\documentclass{amsart}
\newcommand\R{\mathbb R}
\newcommand\N{\mathbb N}
\newcommand\Z{\mathbb Z}
\newcommand\loc{{\text{\upshape loc}}}

\newcommand{\vertiii}[1]{{\left\vert\kern-0.25ex\left\vert\kern-0.25ex\left\vert #1
    \right\vert\kern-0.25ex\right\vert\kern-0.25ex\right\vert}}
\DeclareMathOperator\Tr{Tr}
\newtheorem{thm}{Theorem}

\newtheorem{lem}{Lemma}

\theoremstyle{definition}
\newtheorem{definition}{Definition}
\newtheorem{example}{Example}
\theoremstyle{remark}
\newtheorem{rem}{Remark}

\begin{document}
\title[Nontrivial solutions for the Laplace equation with...]
{Nontrivial solutions for the Laplace equation with a nonlinear Goldstein--Wentzell boundary condition}
\author{Enzo Vitillaro}
\address[E.~Vitillaro]
       {Dipartimento di Matematica e Informatica, Universit\`a di Perugia\\
       Via Vanvitelli,1 06123 Perugia ITALY}
\email{enzo.vitillaro@unipg.it}

\subjclass{35D30, 35J05,35J20,25J25,35J61,35J67}

\keywords{Laplace equation, Laplace--Beltrami operator, existence and multiplicity for nontrivial solutions,Wentzell boundary conditions, Ventcel boundary conditions, Mountain Pass Theorem}


\thanks{
The research was partially supported by the MIUR - PRIN 2022 project ''Advanced theoretical aspects in PDEs and their applications'' (Prot. N.  2022BCFHN2),   by  the INdAM - GNAMPA Project ``Equazioni differenziali alle derivate parziali di tipo misto o dipendenti da campi di vettori'' (Project number CUP\_E53C22001930001), and by ``Progetti Equazioni delle onde con condizioni iperboliche ed acustiche al bordo,  finanziati  con  i Fondi  Ricerca  di Base 2017--2022, della Universit\`a degli Studi di Perugia''.
}

\begin{abstract} The paper deals with the existence and multiplicity of nontrivial solutions for the doubly elliptic problem
$$\begin{cases} \Delta u=0 \qquad &\text{in
$\Omega$,}\\
u=0 &\text{on $\Gamma_0$,}\\
-\Delta_\Gamma u +\partial_\nu u =|u|^{p-2}u\qquad
&\text{on
$\Gamma_1$,}
\end{cases}
$$
where $\Omega$ is a bounded open subset of $\R^N$ ($N\ge
2$) with $C^1$ boundary  $\partial\Omega=\Gamma_0\cup\Gamma_1$, $\Gamma_0\cap\Gamma_1=\emptyset$,
$\Gamma_1$ being nonempty and relatively open on $\Gamma$,  $\mathcal{H}^{N-1}(\Gamma_0)>0$ and $p>2$ being subcritical with respect to Sobolev embedding on $\partial\Omega$.

We prove that the problem admits nontrivial solutions at the  potential--well depth energy level, which is the minimal energy level for nontrivial solutions. We also  prove that the problem has infinitely many solutions at higher energy levels.
\end{abstract}

\maketitle
\section{Introduction and main results} \label{intro}
We deal with the doubly elliptic problem
\begin{equation}\label{1}
\begin{cases} \Delta u=0 \qquad &\text{in
$\Omega$,}\\
u=0 &\text{on $\Gamma_0$,}\\
-\Delta_\Gamma u +\partial_\nu u =|u|^{p-2}u\qquad
&\text{on
$\Gamma_1$,}
\end{cases}
\end{equation}
where $\Omega$ is a bounded open subset of $\R^N$ ($N\ge
2$) with $C^1$ boundary (see \cite{grisvard}). We denote  $\Gamma=\partial\Omega$ and we assume
$\Gamma=\Gamma_0\cup\Gamma_1$, $\Gamma_0\cap\Gamma_1=\emptyset$,
$\Gamma_1$ being nonempty and relatively open on $\Gamma$ (or equivalently $\overline{\Gamma_0}=\Gamma_0$).
Denoting by
$\mathcal{H}^{N-1}$ the  Hausdorff  measure, we assume that   $\mathcal{H}^{N-1}(\overline{\Gamma}_0\cap\overline{\Gamma}_1)=0$ and $\mathcal{H}^{N-1}(\Gamma_0)>0$.
 These
properties of $\Omega$, $\Gamma_0$ and $\Gamma_1$ will be assumed,
without further comments, throughout the paper.
Moreover, in \eqref{1}, we consider $p>2$ and we respectively denote by $\Delta$ and $\Delta_\Gamma$ the Laplace and the Laplace--Beltrami operators, while $\nu$ stands for the outward normal to $\Omega$.

Elliptic equations with nonlinear Neumann boundary conditions, such as problem \eqref{1} without the Laplace--Beltrami term, have a wide literature. Without any aim of completeness, here we refer to \cite{atkinson,MBAHFAS,chfireich,desouza,quitreich,tsungfang}.

Boundary conditions like the one in \eqref{1}, but without the nonlinear source $|u|^{p-2}u$, are known in the literature as generalized Wentzell (sometimes spelled as Vencel) or Goldstein--Wentzell boundary condition, since they have been subject of several papers in the framework of linear evolutions problems. See for example \cite{grecoviglialoro,quarteroni,nicaisemazzucato,Rom2014DCDS,vazvitHLB} and \cite{lionstata}, to which we refer for  the physical motivations of this kind of problems.

On the other hand, to the author's knowledge, a Goldtsein--Wentzell boundary condition with a nonlinear source like $|u|^{p-2}u$ in connection with the Laplace equation has never been considered in the literature. The motivation for studying it comes from a series of papers by the author concerning the wave equation with hyperbolic dynamical boundary conditions with boundary damping and source terms. The prototype of this kind of problems is the evolutionary boundary value problem
\begin{equation}\label{2}
\begin{cases} u_{tt}-\Delta u=0 \qquad &\text{in
$(0,\infty)\times\Omega$,}\\
u=0 &\text{on $(0,\infty)\times \Gamma_0$,}\\
u_{tt}+\partial_\nu u-\Delta_\Gamma
u+|u_t|^{m-2}u_t=|u|^{p-2}u\qquad &\text{on
$(0,\infty)\times \Gamma_1$,}
\end{cases}
\end{equation}
where $u=u(t,x)$, $t\ge 0$, $x\in\Omega$, $\Delta=\Delta_x$ denotes
the Laplacian operator respect to the space variable.
 Its associated initial--value problem was introduced in \cite{AMS} and then studied, as a particular case, in \cite{Dresda1, Dresda2} and \cite{Dresda3}. We refer to \cite{goldsteingisele,Dresda1} for the physical derivation of the problem, describing the vibrations of a membrane with a part of the boundary carrying a linear density of  kinetic energy.

 In order to give clear--cut criteria on the initial data to discriminate between global existence and blow--up for solutions of \eqref{2} it is useful to know if it possesses nontrivial stationary solutions, which turns out to be solutions of \eqref{1}, at some specific energy level.

 In particular in the present paper we shall consider the case when the nonlinearity $|u|^{p-2}u$ is sub--critical with respect to the
 Sobolev Embedding $H^1(\Gamma)\hookrightarrow L^p(\Gamma)$, that is we shall assume that
 \begin{equation}\label{3}
 2<p<r,\qquad\text{where}\quad r=
 \begin{cases}
 \frac{2(N-1)}{N-3} \quad&\text{if $N\ge 4$,}\\
 \quad \infty  \quad&\text{if $N=2,3$.}\\
 \end{cases}
 \end{equation}
Moreover, when dealing with problem \eqref{3}, we shall also assume that
\begin{equation}\label{15}
 m>1,\qquad p\le 1+r/\overline{m}',\quad \text{where}\quad \overline{m}:=\max\{2,m\},
 \end{equation}
 the last assumption being related with well--posedness issues, see the papers quoted above.
 \begin{footnote}{By the way assumption \eqref{15} may be skipped when dealing with stationary solutions, but we prefer to keep it  to avoid re--discussing here problem~\eqref{2}}\end{footnote}
  We also remark that, although also the case $p\ge r$ (when $N\ge 4$) was considered there, only the case $p<r$ is of interest when dealing with the dichotomy between global existence and blow--up, see \cite[Remark~1, p.6]{Dresda3}.

 To state our main results we first introduce some basic notation. In the sequel we shall
identify $L^p(\Gamma_1)$, for $1\le p\le \infty$, with its isometric image in $L^p(\Gamma)$,
that is
\begin{equation}\label{4}
L^p(\Gamma_1)=\{u\in L^p(\Gamma):\,\, u=0\,\,\text{a.e. on
}\,\,\Gamma_0\}.
\end{equation}
Moreover we shall denote by $\Tr$ the trace operator from $H^1(\Omega)$ onto $H^{1/2}(\Gamma)$ and, for simplicity of notation,
$\Tr u=u_{|\Gamma}$.

We introduce the Hilbert spaces $H^0 = L^2(\Omega)\times L^2(\Gamma_1)$ and
\begin{equation}\label{5}
H^1 = \left\{(u,v)\in H^1(\Omega)\times H^1(\Gamma): v=u_{|\Gamma}, v=0
\,\,\ \text{on $\Gamma_0$}\right\},
\end{equation}
with the topologies inherited from the products. For the sake of
simplicity we shall identify, when useful, $H^1$ with its isomorphic
counterpart
\begin{equation}\label{5bis}
H^1_{\Gamma_0}(\Omega,\Gamma)=\{u\in H^1(\Omega): u_{|\Gamma}\in H^1(\Gamma)\cap
L^2(\Gamma_1)\},
\end{equation}
studied for example in \cite{pucvit}, through the identification $(u,u_{|\Gamma})\mapsto
u$. So we shall write, without further mention, $u\in H^1$ for
functions defined on $\Omega$. Moreover we shall drop the notation
$u_{|\Gamma}$, when useful, so we shall write $\|u\|_{L^2(\Gamma)}$
and so on, for $u\in H^1$, referring to the restriction of the Hausdorff measure $\mathcal{H}^{N-1}$ to measurable subsets of $\Gamma$.
We shall also drop the notation  $d\mathcal{H}^{N-1}$ in boundary integrals, so writing $\int_\Gamma u=\int_\Gamma u\,d\mathcal{H}^{N-1}$.

By assumption \eqref{3} we can introduce in $H^1$ the nonlinear functional $I\in C^1(H^1)=C^1(H^1;\R)$ defined by
\begin{footnote}{here $\nabla_\Gamma$
denotes the Riemannian gradient on $\Gamma$ and $|\cdot|_\Gamma$,
the norm associated to the Riemannian scalar product on the tangent
bundle of $\Gamma$. See Section \ref{section 2}.}\end{footnote}

\begin{equation}\label{6}
 I(u)=\tfrac 12 \int_\Omega |\nabla u|^2 +\tfrac 12
\int_{\Gamma_1} |\nabla_\Gamma  u|_\Gamma^2-\tfrac 1p \int_{\Gamma_1}|u|^p,
\end{equation}
which represents the potential energy associated to problem \eqref{2}. For this reason we shall call it the \emph{energy functional} when dealing with  \eqref{1}.

We also introduce the potential--well depth $d$ given by
\begin{equation}\label{7}
 d=\inf_{u\in H^1, u_{|\Gamma}\not\equiv 0}\sup_{\lambda>0} I(\lambda u)=\inf_{u\in H^1\setminus\{0\}}\sup_{\lambda>0} I(\lambda u),
\end{equation}
noticing that the identity between the two infima in \eqref{7} is essentially trivial  and that we shall prove  that $d>0$.

Our first main result shows that problem \eqref{1} admits  nontrivial weak solutions, see Definition~\ref{Definition 3} below, coinciding with critical points of the functional $I$, at the positive energy level $d$. We shall also recognize that they are stationary weak solutions of \eqref{3} provided this class of solutions is well--defined, see Definition~\ref{Definition 2} below.

\begin{thm}\label{Theorem 1} When \eqref{3} holds problem \eqref{1} has at least a couple $(u,-u)$ of antipodal weak solutions in $H^1$ such that $I(u)=I(-u)=d>0$. When \eqref{15} holds  they are also stationary weak solutions of problem \eqref{2}.

Moreover $d$ coincides with the Mountain Pass level of the functional $I$, that is
$$d=\inf_{\sigma\in\Sigma} \max_{t\in [0,1]}I(\sigma(t)), \quad\text{where}\quad \Sigma=\{\sigma\in C([0,1];H^1):\,\sigma(0)=0, I(\sigma(1))<0\}.
$$
\end{thm}
Theorem~\ref{Theorem 1} will be proved by applying  a  variant, maybe less well--known than other ones, of the Mountain Pass Theorem, explicitly given in \S~\ref{section 2}.

To show the relevance of the potential--well depth $d$, beside its interest in evolution problems, we give some relevant properties of weak solutions of \eqref{1} having energy $d$. At first we introduce the norm $B$ of the bounded linear trace operator from $H^1$ to $L^p(\Gamma_1)$, that is
\begin{equation}\label{7bis}
  B=\sup_{u\in H^1\setminus\{0\}}\frac {\|u\|_{L^p(\Gamma_1)}}{\left(\|\nabla u\|_{L^2(\Omega)}^2+\|\nabla_\Gamma u\|_{L^2(\Gamma_1)}^2\right)^{1/2}},
\end{equation}
noticing that we shall prove that $B<\infty$. We can then state our second main result.
\begin{thm}\label{Theorem 1.1} Let \eqref{3} hold and set
\begin{equation}\label{7ter}
 \lambda_1=B^{-p/(p-2)},\qquad\text{and}\quad  \lambda_2=B^{-2/(p-2)}.
\end{equation}
Then we have
\begin{equation}\label{43}
 d=\left(\tfrac 12 -\tfrac 1p\right)\lambda_1^2=\left(\tfrac 12 -\tfrac 1p\right)\lambda_2^p.
\end{equation}
Moreover, if $u$ is a weak solution of \eqref{1} with $I(u)=d$ we also have
\begin{equation}\label{7quater}
\|\nabla u\|_{L^2(\Omega)}^2+\|\nabla_\Gamma u\|_{L^2(\Gamma_1)}^2=\lambda_1^2\qquad\text{and}\quad \|u\|_{L^p(\Gamma_1)}=\lambda_2.
\end{equation}
Finally weak solutions at the energy level $d$ are least energy non--trivial solutions of \eqref{1}, that is for any non--trivial weak solutions $u$ of \eqref{1} one has $I(u)\ge d$.
\end{thm}
The proof of  the minimality of the energy of solutions  at level $d$, stated in  Theorem~\ref{Theorem 1.1}, is of  elementary nature. Hence it is  easier than most proofs in the literature for internal sources, see for example \cite{jeanjeantanaka2,jeanjeantanaka1,tanakaMP}. By the way the homogeneity of the source $|u|^{p-2}u$ allows this simple approach.

Finally, to show that the minimality asserted in Theorem~\ref{Theorem 1.1} is of some use, since there are solutions at an higher level, we give our last main result, which can be also of independent interest.
\begin{thm}\label{Theorem 2} When \eqref{3} holds there is a sequence $(u_n)_n$ of nontrivial weak solutions of \eqref{1} such that $I(u_n)\to\infty$ as $n\to\infty$.
\end{thm}
The proof of Theorem~\ref{Theorem 2} relies on applying the $\Z_2$--version of the Mountain Pass Theorem in a different variational setting, which turns out to be equivalent to the one illustrated  in this Section. See \S~\ref{section 4} for details.

We would like to mention that, although in the paper we give Theorems~\ref{Theorem 1}--\ref{Theorem 2} for the prototype nonlinearity $f(x,u)=|u|^{p-2}u$, they can be easily extended to the  problem
$$\begin{cases} \Delta u=0 \qquad &\text{in
$\Omega$,}\\
u=0 &\text{on $\Gamma_0$,}\\
-\Delta_\Gamma u +\partial_\nu u =f(x,u)\qquad
&\text{on
$\Gamma_1$,}
\end{cases}
$$
under suitable  assumptions on $f$ which have been widely used in the literature. Here, for the sake of simplicity, we preferred to concentrate on the prototype problem.

Moreover clearly Theorems~\ref{Theorem 1}--\ref{Theorem 2} also show the existence of  stationary solutions for other evolution problems not considered in the present paper, which can be of future interest.

The paper is organized as follows: in Section~\ref{section 2} we shall give all preliminaries needed in the paper.  Section~\ref{section 3} will be devoted to prove Theorems~\ref{Theorem 1} and \ref{Theorem 1.1}, while in Section~\ref{section 4} we shall prove Theorem~\ref{Theorem 2}.

\section{Preliminaries} \label{section 2}
\subsection{Notation.}\label{section 2.1}
We shall adopt the standard notation for
(real) Lebesgue and Sobolev spaces in $\Omega$, referring to
\cite{adams}. For simplicity we shall denote by $\|\cdot\|_{\tau}$, for $1\le \tau\le \infty$, the norms in $L^\tau(\Omega)$ and in
$L^\tau(\Omega;\R^N)$.

Given a Banach space
$X$ we shall denote by  $X'$  its dual and by $\langle
\cdot,\cdot\rangle_X$ the duality product between them. Moreover we shall use the standard notation for $X$--valued Lebesgue and Sobolev spaces in a real interval. When another Banach space $Y$ is given we shall denote by
$\mathcal{L}(X,Y)$  the space of bounded linear operators between $X$ and $Y$ and by $\|\cdot\|_{\mathcal{L}(X,Y)}$ the standard norm on it.

\subsection{Function spaces and Riemannian operators on $\Gamma$.} \label{section 2.2}
Lebesgue spaces on $\Gamma$ and $\Gamma_1$ will be intended with respect to (the restriction to measurable subset of them of) the Hausdorff measure
$\mathcal{H}^{N-1}$, and for simplicity we shall denote, for $1\le \tau\le \infty$,
$\|\cdot\|_{\tau,\Gamma}=\|\cdot\|_{L^\tau(\Gamma)}$ and $\|\cdot\|_{\tau,\Gamma_1}=\|\cdot\|_{L^\tau(\Gamma_1)}$.

 Sobolev spaces on $\Gamma$ and on its relatively open subsets are classical objects, and  we shall use the standard notation for them. We refer to \cite{grisvard} for their definition in the present case in which $\Gamma$ is merely $C^1$.

  Since $\Gamma$ is $C^1$  it inherits from $\R^N$ the structure of a Riemannian $C^1$ manifold, see \cite{sternberg}, so in the sequel we shall use some notation of geometric nature, which is quite common when $\Gamma$ is smooth,  see \cite{Boothby, hebey, jost, taylor}, and which can be easily extended to the $C^1$ case, see for example \cite{mugnvit}.  Moreover, since $\Gamma_1$ is relatively open on $\Gamma$, this notation will apply (by restriction) to it, without further mention.

We shall denote by $T(\Gamma)$ and $T^*(\Gamma)$ the tangent and cotangent bundles, and  by $(\cdot,\cdot)_\Gamma$ the Riemannian metric inherited from $\R^N$, given in local coordinates by $(u,v)_\Gamma=g_{ij}u^i v^j$ for all $u,v\in T(\Gamma)$
(here and in the sequel the summation convention being in use). The metric induces the fiber--wise defined musical isomorphisms $\flat:T(\Gamma)\to T^*(\Gamma)$ and $\sharp=\flat^{-1}:T^*(\Gamma)\to T(\Gamma)$ defined by
$\langle \flat u,v\rangle_{T(\Gamma)}=(v,u)_\Gamma$ for $u,v\in T(\Gamma)$,
where $\langle \cdot,\cdot\rangle_{T(\Gamma)}$ denotes the fiber-wise defined duality pairing.
The induced bundle metric on $T^*(\Gamma)$, still denoted by $(\cdot,\cdot)_\Gamma$, is then defined  by
the formula $(\alpha,\beta)_\Gamma=\langle \alpha,\sharp\beta\rangle_{T(\Gamma)}$ for all $\alpha,\beta\in T^*(\Gamma)$, so that
\begin{equation}\label{+}
  (\alpha,\beta)_\Gamma=(\sharp \beta, \sharp\alpha)_\Gamma,\qquad \text{for all $\alpha,\beta\in T^*(\Gamma)$}.
\end{equation}
By $|\cdot|_\Gamma^2= (\cdot,\cdot)_\Gamma$  we shall denote the associated bundle norms on $T(\Gamma)$ and $T^*(\Gamma)$.

Denoting by $d_\Gamma$ the standard differential on $\Gamma$, the Riemannian gradient operator $\nabla_\Gamma$ is defined  by
setting, for $u\in C^1(\Gamma)$ and thus by density for $u\in H^1(\Gamma)$,
$\nabla_\Gamma u=\sharp d_\Gamma u$,
so $\nabla_\Gamma u=g^{ij}\partial_ju\partial_i$
 in local coordinates, where $(g^{ij})=(g_{ij})^{-1}$. By \eqref{+} one trivially gets that
$(\nabla_\Gamma u,\nabla_\Gamma v)_\Gamma=(d_\Gamma u,d_\Gamma v)_\Gamma$ for all $u,v\in H^1(\Gamma)$, so in the sequel the use of vectors or forms is optional.

It is well known, see for example \cite[Chapter 3]{mugnvit}, that $H^1(\Gamma)$ can be equipped with the equivalent norm $\|\cdot\|_{H^1(\Gamma)}$ given by
$$\|u\|_{H^1(\Gamma)}^2=\|u\|_{2,\Gamma}^2+\|\nabla_\Gamma u\|_{2,\Gamma}^2,\quad\text{where}\quad \|\nabla_\Gamma u\|_{2,\Gamma}^2:=\int_\Gamma |\nabla_\Gamma u|_\Gamma^2.$$
In the sequel we shall also deal with the closed subspace of $H^1(\Gamma)$
\begin{equation}\label{++}
H^1_{\Gamma_0}(\Gamma)=\{u\in H^1(\Gamma): u=0\quad\text{a.e. on $\Gamma_0$}\}
\end{equation}
endowed with the  norm $\|\cdot\|_{H^1(\Gamma)}$, which is then an Hilbert space. Since for all $u\in H^1_{\Gamma_0}(\Gamma)$ one has
$\nabla_\Gamma u=0$ a.e. on $\Gamma\setminus\overline{\Gamma_1}$ and $\mathcal{H}^{N-1}(\overline{\Gamma_0}\cap\overline{\Gamma_1})=0$, we have
\begin{equation}\label{+++}
\|u\|_{H^1(\Gamma)}^2=\|u\|_{2,\Gamma_1}^2+\|\nabla_\Gamma u\|_{2,\Gamma_1}^2\qquad\text{for all $u\in H^1_{\Gamma_0}(\Gamma)$,}
\end{equation}
where $\|\nabla_\Gamma u\|_{2,\Gamma_1}^2:=\int_{\Gamma_1} |\nabla_\Gamma u|_\Gamma^2$.

\begin{rem}\label{rem*} Although the definition of the space $H^1_{\Gamma_0}(\Gamma)$ given above is adequate for our purpose, we would like to point out two characterizations of it in two different geometrical situations:
\renewcommand{\labelenumi}{{\roman{enumi})}}
\begin{enumerate}
\item when $\overline{\Gamma_0}\cap\overline{\Gamma_1}=\emptyset$, so both $\Gamma_0$ and $\Gamma_1$ are relatively open, by identifying the elements of $H^1(\Gamma_i)$, $i=0,1$, with their trivial extensions to $\Gamma$, one easily gets the splitting $H^1(\Gamma)=H^1(\Gamma_0)\oplus H^1(\Gamma_1)$, and consequently  $H^1_{\Gamma_0}(\Gamma)$ is isometrically isomorphic to $H^1(\Gamma_1)$;
\item  when $\overline{\Gamma_0}\cap\overline{\Gamma_1}\not=\emptyset$ such a characterization is false, since one easily sees that the characteristic function $\chi_{\Gamma_1}$ of $\Gamma_1$ does not belong to $H^1_{\Gamma_0}(\Gamma)$, while its restriction to $\Gamma_1$ trivially belongs to $H^1(\Gamma_1)$. Indeed in this case the elements of $H^1_{\Gamma_0}(\Gamma)$ ''vanish'' at the relative boundary $\partial\Gamma_1=\overline{\Gamma_0}\cap\overline{\Gamma_1}$ of $\Gamma_1$ on $\Gamma$, although such a notion can be made more precise only when $\partial\Gamma_1$ is regular enough. For example, when $\Gamma$ is smooth and $\overline{\Gamma_1}$ is a manifold with boundary $\partial\Gamma_1$,   see \cite[\S 5.1]{taylor}, $H^1_{\Gamma_0}(\Gamma)$ is isometrically isomorphic to  the space
    $$H^1_0(\Gamma_1):=\overline{C^\infty_c(\Gamma_1)}^{\|\cdot\|_{H^1(\Gamma_1)}}.$$
\end{enumerate}
\end{rem}
The Laplace--Beltrami operator $\Delta_\Gamma$ can be defined in a geometrically elegant way by using $\nabla_\Gamma$ and the Riemannian divergence operator, as in \cite[\S~2.3]{mugnvit}, at least when $\Gamma$ is $C^2$. To avoid the need of introducting Sobolev spaces of tensor fields we shall adopt here a less elegant approach.
Indeed we set, when $\Gamma$ is $C^2$ and $u\in C^2(\Gamma')$, $\Gamma'\subset\Gamma$ relatively open,
\begin{equation}\label{runa1}
\Delta_{\Gamma} u=g^{-1/2}\partial_i( g^{1/2}g^{ij}\partial_j u), \quad\text{where $g=\det (g_{ij})$,}
\end{equation}
in local coordinates. Since $g$, $g^{ij}$ are continous and $\Gamma$ is compact, formula \eqref{runa1} extends by density to $u\in H^2(\Gamma)$, so defining an operator $-\Delta_\Gamma\in\mathcal{L}(H^2(\Gamma);L^2(\Gamma))$, which restricts to $-\Delta_\Gamma\in\mathcal{L}(H^2(\Gamma');L^2(\Gamma'))$ for relatively open subsets $\Gamma'$ of $\Gamma$. Since $\Gamma$ is compact, by \eqref{runa1}, integrating by parts and using a $C^2$ partition of the unity one gets that
\begin{equation}\label{runa2}
-\int_\Gamma \Delta_\Gamma u v=\int_\Gamma (\nabla_\Gamma u,\nabla_\Gamma v)_\Gamma\quad\text{for all $u\in H^2(\Gamma)$ and $v\in H^1(\Gamma)$.}
\end{equation}
Formula \eqref{runa2} motivates the definition of the operator $-\Delta_\Gamma\in\mathcal{L}(H^1(\Gamma);H^{-1}(\Gamma))$, also when $\Gamma$ is merely $C^1$, given by
\begin{equation}\label{runa3}
\langle -\Delta_\Gamma u,v\rangle _{H^1(\Gamma)}=\int_\Gamma (\nabla_\Gamma u, \nabla_\Gamma v)_\Gamma\quad\text{for all $u,v\in H^1(\Gamma)$.}
\end{equation}
By density, when $\Gamma$ is $C^2$, the so defined operator is the unique  extension of $-\Delta_\Gamma\in\mathcal{L}(H^2(\Gamma);L^2(\Gamma))$.

In \S~\ref{section 4} we shall deal with the realization of $-\Delta_\Gamma$ between the space $H^1_{\Gamma_0}(\Gamma)$ and its dual. The different nature of the space $H^1_{\Gamma_0}(\Gamma)$  in the two cases i) and ii) has been pointed out in Remark~\ref{rem*}. To explain the definition of the realization we shall give we recall that, when $\overline{\Gamma_0}\cap\overline{\Gamma_1}=\emptyset$, so $\Gamma_1$ is compact, formula \eqref{runa2} holds, when $\Gamma$ is $C^2$, also when replacing $\Gamma$ with $\Gamma_1$, so making natural to set $-\Delta_{\Gamma_1}\in\mathcal{L}(H^1(\Gamma_1);H^{-1}(\Gamma_1))$ by
\begin{equation}\label{runa4}
\langle -\Delta_{\Gamma_1} u,v\rangle _{H^1(\Gamma_1)}=\int_{\Gamma_1} (\nabla_\Gamma u, \nabla_\Gamma v)_\Gamma\quad\text{for all $u,v\in H^1(\Gamma_1)$.}
\end{equation}
When $\overline{\Gamma_0}\cap\overline{\Gamma_1}\not=\emptyset$, $\Gamma$ is smooth and $\overline{\Gamma_1}$ is a manifold with boundary $\partial\Gamma_1$, formula \eqref{runa2} fails to hold on $\Gamma_1$, since a boundary integral on $\partial\Gamma_1$ appears. On the other hand, taking into account the homogeneous Dirichlet boundary condition in the space $H^1_0(\Gamma_1)$, it is natural to set $-\Delta_{\Gamma_{1D}}\in\mathcal{L}(H^1_0(\Gamma_1);H^{-1}(\Gamma_1))$ by
\begin{equation}\label{runa5}
\langle -\Delta_{\Gamma_{1D}} u,v\rangle _{H^1_0(\Gamma_1)}=\int_{\Gamma_1} (\nabla_\Gamma u, \nabla_\Gamma v)_\Gamma\quad\text{for all $u,v\in H^1_0(\Gamma_1)$.}
\end{equation}
Hence, taking into account the characterizations given in Remark~\ref{rem*}, to simultaneously deal with the two cases i) and ii), in the sequel we shall deal with the operator $-\Delta_{\Gamma_{1(D)}}\in\mathcal{L}(H^1_{\Gamma_0}(\Gamma);[H^1_{\Gamma_0}(\Gamma)]')$
defined by
\begin{equation}\label{**}
\langle -\Delta_{\Gamma_{1(D)}} u,v\rangle _{H^1_{\Gamma_0}(\Gamma)}=\int_{\Gamma_1} (\nabla_\Gamma u, \nabla_\Gamma v)_\Gamma\quad\text{for all $u,v\in H^1_{\Gamma_0}(\Gamma_1)$,}
\end{equation}
noticing that, by \eqref{runa3}, $-\Delta_{\Gamma_1(D)}u=-\Delta_{\Gamma}u_{|H^1_{\Gamma_0}(\Gamma)}$ for all $u\in H^1_{\Gamma_0}(\Gamma)$.
\subsection{The space $H^1$.} \label{section 2.3} We recall, see \cite[Lemma 1, p. 2147]{vazvitHLB} which trivially extends to $\Gamma$ of class $C^1$, that the space
$$H^1(\Omega;\Gamma)=\{(u,v)\in H^1(\Omega)\times H^1(\Gamma): v=u_{|\Gamma}\},$$
with the topology inherited from the product, can be identified with the space $\{u\in H^1(\Omega): u_{|\Gamma}\in H^1(\Gamma)\}$ and equivalently equipped with the norm $\|\cdot\|_{H^1(\Omega,\Gamma)}$ given by
$$\|u\|^2_{H^1(\Omega,\Gamma)}=\|\nabla u\|_2^2+\|\nabla_\Gamma u\|_{2,\Gamma}^2+\|u\|_{2,\Gamma}^2.$$
The identification made in \S~\ref{intro} between the spaces $H^1$ and $H^1_{\Gamma_0}(\Omega,\Gamma)$, respectively defined  by \eqref{5} and \eqref{5bis}, is a simple consequence of the identification above, and, by \eqref{+++},  $H^1$ can be equivalently equipped with the norm $\vertiii\cdot_{H^1}$ given by
\begin{equation}\label{8}
 \vertiii u_{H^1}^2= \|\nabla u\|_2^2+\|\nabla_\Gamma u\|_{2,\Gamma_1}^2+\|u\|_{2,\Gamma_1}^2.
\end{equation}
On the other hand, to get advantage of the assumption $\mathcal{H}^{N-1}(\Gamma_0)>0$, made in the present paper, we point out the following well--known result, the proof of which  is given only for the reader's convenience.
\begin{lem}\label{lemma 0} Let $\mathcal{H}^{N-1}(\Gamma_0)>0$. Then, setting, for $u,v\in H^1$,
\begin{equation}\label{10}
  (u,v)_{H^1}=\int_\Omega\nabla u\nabla v+\int_{\Gamma_1}(\nabla_\Gamma u,\nabla_\Gamma v)_\Gamma \quad\text{and}\quad \|\cdot\|_{H^1}=(\cdot,\cdot)_{H^1}^{1/2},
\end{equation}
 $\|\cdot\|_{H^1}$ defines on $H^1$ a norm equivalent to  $\vertiii\cdot_{H^1}$.
\end{lem}
\begin{proof}By combining \cite[Chapter 2, Theorem 2.6.16, p. 75]{ziemer} and \cite[Chapter 4, Corollary 4.5.2, p. 195]{ziemer} one gets the following Poincar\'{e}--type inequality: there is a positive constant $c_1=c_1(\Omega,\Gamma_0)$ such that
\begin{equation}\label{9}
  \|u\|_2\le c_1 \|\nabla u\|_2\qquad\text{for all $u\in H^1$.}
\end{equation}
By the Trace Theorem there is a positive constant $c_2=c_2(\Omega,\Gamma_0)$ such that
\begin{equation}\label{11}
  \|u\|_{2,\Gamma_1}\le c_2 \| u\|_{H^1(\Omega)}\qquad\text{for all $u\in H^1(\Omega)$,}
\end{equation}
where $\|\cdot\|_{H^1(\Omega)}$ is the standard norm of $H^1(\Omega)$. Since $H^1\subset H^1(\Omega)$, by combining \eqref{8}, \eqref{9} and \eqref{11} we get
$$\vertiii u_{H^1}^2\le \|\nabla u\|_2^2+\|\nabla_\Gamma u\|_{2,\Gamma_1}^2+ c_2(1+c_1)\|\nabla u\|_2^2\le c_3\|u\|_{H^1}^2$$
for all $u\in H^1$, where $c_3=1+c_2(1+c_1)$, from which the statement trivially follows.
\end{proof}
\subsection{Some results from Critical Point Theory} \label{section 2.4} We now recall some well--known notions of Critical Point Theory for a functional $\mathcal{I}\in C^1(X)=C^1(X;\R)$ on any Banach space $X$ with norm $\|\cdot\|_X$. By $\mathcal{I}'\in C(X;X')$ we shall denote the Fr\'{e}chet differential of $\mathcal{I}$ and (PS) will stands, in short,  for Palais--Smale. See \cite{rabinowitz}.
\begin{definition}\label{Definition 1} Let $\mathcal{I}\in C^1(X)$. We say that a sequence $(u_n)_n$ in $X$  is a  (PS) sequence if $(\mathcal{I}(u_n))_n$ is bounded and $\mathcal{I}'(u_n)\to 0$ in $X'$. We also say that $\mathcal{I}\in C^1(X)$ satisfies the (PS) condition if any (PS) sequence has a (strongly) convergent subsequence.
\end{definition}
The following result is nothing but a well--known version of  the celebrated Mountain Pass Theorem, see \cite[Chapter 1, p. 4]{rabinowitz}.
\begin{thm}\label{Theorem 3} Let $\mathcal{I}\in C^1(X)$ satisfies the (PS) condition and
\renewcommand{\labelenumi}{{\roman{enumi})}}
\begin{enumerate}
\item $\mathcal{I}(0)=0$;
\item there are $\rho,\alpha>0$ such that $\mathcal{I}(u)\ge \alpha$ for all $u\in X$ such that $\|u\|_X=\rho$;
\item there is $l\in X$ such that $\|l\|_X>\rho$ and $\mathcal{I}(l)\leq 0$.
\end{enumerate}
Then $\mathcal{I}$ possesses a critical value $c_l\ge \alpha$ given by
\begin{equation}\label{13}
c_l=\inf_{\sigma\in \Sigma_l}\max_{t\in [0,1]}\mathcal{I}(\sigma(t)), \text{where } \Sigma_l=\{\sigma\in C([0,1];X): \sigma (0)=0, \sigma(1)=l\}.
\end{equation}
\end{thm}
It is rarely pointed out  in textbooks that the critical level $c_l$ above may depend on $l$, as the following trivial example shows.
\begin{example} Let $X=\R$ and
$$\mathcal{I}(x)=
\begin{cases}
x^2-x^4 &\quad\text{if $x\ge 0$,}\\
x^2-2x^4&\quad\text{if $x<0$.}
\end{cases}$$
Trivially $\mathcal{I}\in C^1(\R)$ satifies the (PS) condition as well as assumptions i)--iii). Moreover its critical points are exactly $x=0,-1/2,\sqrt 2/2$ from which one easily sees that
$$c_l=
\begin{cases}
\mathcal{I}(\sqrt 2/2)=1/4 &\quad\text{if $l>0$,}\\
\mathcal{I}(-1/2)=3/16 &\quad\text{if $l<0$.}
\end{cases}$$
\end{example}
Since in the present paper we are interested in characterizing our critical level as the potential--well depth of the functional $I$, we now state a less known variant of the Mountain Pass Theorem under slighty more restrictive assumptions on the functional which look
similar (although not identical) to the assumptions of \cite[Theorem 2.1, p. 354]{AmbrosettiRabinowitz}, this one being the first version of this celebrated result.
\begin{thm}\label{Theorem 4} Let $\mathcal{I}\in C^1(X)$ satisfies the (PS) condition, assumptions i)--iii) in Theorem~\ref{Theorem 3} and and
\renewcommand{\labelenumi}{{\roman{enumi})}}
\begin{itemize}
\item[iv)] $\mathcal{I}(u)>0$ for all $u\in X\setminus\{0\}$ such that $\|u\|_X\le\rho$.
\end{itemize}
Then $\mathcal{I}$ possesses a critical value $c\ge \alpha$ given by
\begin{equation}\label{14}
c=\inf_{\sigma\in \Sigma}\max_{t\in [0,1]}\mathcal{I}(\sigma(t)), \text{where } \Sigma=\{\sigma\in C([0,1];X): \sigma (0)=0,\, \mathcal{I}(\sigma(1))<0\}.
\end{equation}
\end{thm}
\begin{proof} One can  repeat the proof of \cite[Theorem 2.1, p. 354]{AmbrosettiRabinowitz} verbatim, since by assumption iv) for any  $\sigma\in\Sigma$ one has $\|\sigma(1)\|_X>\rho$. Alternatively one can also deduce the statement from Theorem~\ref{Theorem 3} by the following simple argument. By assumption iv) one has $\Sigma=\bigcup_{l\in A}\Sigma_l$, where $A=\{l\in X: \|l\|_X>\rho\,\,\text{and}\,\,\mathcal{I}(l)<0\}$. Hence, by \eqref{13} and \eqref{14} we have $c=\inf_{l\in A}c_l$. Hence, since by assumption ii) one has $\alpha\le c<\infty$, there is a sequence $(l_n)_n$ in $A$ such that $c_{l_n}\to c$. By Theorem~\ref{Theorem 3} there ia a corresponding sequence $(u_n)_n$ in $X$ such that $\mathcal{I}(u_n)=c_{l_n}$ and $\mathcal{I}'(u_n)=0$, which is then a (PS) sequence and consequently, up to a subsequence, $u_n\to u$. Then $\mathcal{I}'(u)=0$ and $\mathcal{I}(u)=c$, concluding the proof.
\end{proof}
In the sequel we shall also use the following well--known $\Z_2$--version of the Mountain Pass Theorem, see \cite[Chapter 9, Theorem 9.12, p. 55  and Proposition 9.33, p. 58]{rabinowitz}.
\begin{thm}\label{Theorem 5} Let $X$ be  infinite dimensional and $\mathcal{I}\in C^1(X)$ be even, satisfying the (PS) condition, assumptions i)--ii) of Theorem~\ref{Theorem 3} and
\begin{itemize}
\item[v)] for each finite dimensional subspace $Y$ of $X$ there is $R_Y>0$ such that $\mathcal{I}(u)\le 0$ for all $u\in Y$ such that $\|u\|_X>R_Y$.
\end{itemize}
Then $\mathcal{I}$ possesses a sequence $(u_n)_n$ of critical points such that $\mathcal{I}(u_n)\to\infty$.
\end{thm}
\section{Existence of solutions of \eqref{1} at level $d$ and their minimality.} \label{section 3}
The aim of this section is to prove Theorems~\ref{Theorem 1} and \ref{Theorem 1.1}.
We start by recalling  what we mean by a weak solution of \eqref{2}, referring to \cite[\S 2.2 and Definition 3.1, p. 4896]{Dresda2}. We shall also make  precise the use of the term "stationary" when referring to them.
\begin{definition}\label{Definition 2} Let \eqref{3} and \eqref{15} hold.  A \emph{weak solution  of \eqref{2}} is
\begin{equation}\label{16}
  u\in L^\infty_\loc([0,\infty);H^1)\cap W^{1,\infty}_\loc ([0,\infty);H^0),\quad (u_{|\Gamma_1})_t\in L^m_\loc (0,\infty); L^m(\Gamma_1))
\end{equation}
satisfying the distribution identity
\begin{multline}\label{17}
  \int_0^\infty\left[-\int_\Omega u_t\psi_t-\int_{\Gamma_1}(u_{|\Gamma_1})_t(\psi_{|\Gamma_1})_t \right.\\
  \left. +\int_\Omega \nabla u \nabla\psi +\int_{\Gamma_1}(\nabla_\Gamma u,\nabla_\Gamma\psi)_\Gamma -\int_{\Gamma_1}|u|^{p-2}u\psi \right]=0
\end{multline}
for all $\psi\in C_c((0,\infty);H^1)\cap C^1_c((0,\infty); H^0)$ such that $(\psi_{|\Gamma_1})_t\in L^m_\loc ((0,\infty); L^m(\Gamma_1))$.
We say that $u$ is \emph{stationary} if  $u(t)\equiv u_0\in H^1$ in $(0,\infty)$.
\end{definition}
We also make precise what we mean by weak solutions of \eqref{1}.
\begin{definition}\label{Definition 3} Let $2<p<r$. A \emph{weak solution} of \eqref{1} is $u\in H^1$ such that
\begin{equation}\label{18}
  \int_\Omega \nabla u\nabla \phi+\int_{\Gamma_1}(\nabla_\Gamma u,\nabla_\Gamma \phi)_\Gamma-\int_{\Gamma_1}|u|^{p-2}u\phi=0\quad\text{for all $\phi\in H^1$.}
\end{equation}
\end{definition}
Actually weak solutions of \eqref{1} and stationary weak solutions of \eqref{2} coincide when they are both defined, as the following result shows.
\begin{lem}\label{Lemma 1}
Let  \eqref{3}, \eqref{15} hold and $u_0\in H^1$. Then $u\equiv u_0$ is a stationary weak  solution of \eqref{2} if and only if $u_0$ is a weak solution of \eqref{1}.
\end{lem}
\begin{proof}If $u_0$ is a weak solution of \eqref{1} by \eqref{18} one immediately gets that $u\equiv u_0$ satisfies \eqref{17}. To prove the converse we recall that, by \cite[Lemma 3.3, p. 4896]{Dresda2}, any weak solution $u$ of \eqref{2} satisfies the alternative form of the  distribution identity \eqref{17}
\begin{multline}\label{19}
  \left[\int_\Omega u_t\psi+\int_{\Gamma_1}(u_{|\Gamma_1})_t\psi\right]_0^T+\int_0^T\left[-\int_\Omega u_t\psi_t-\int_{\Gamma_1}(u_{|\Gamma_1})_t(\psi_{|\Gamma_1})_t \right.\\
  \left. +\int_\Omega \nabla u \nabla\psi +\int_{\Gamma_1}(\nabla_\Gamma u,\nabla_\Gamma\psi)_\Gamma -\int_{\Gamma_1}|u|^{p-2}u\psi \right]=0
\end{multline}
for all $T>0$ and $\psi\in C([0,T];H^1)\cap C^1([0,T]; H^0)$, $(\psi_{|\Gamma_1})_t\in L^m ((0,T)\times \Gamma_1$.
Hence, when $u\equiv u_0\in H^1$ is a stationary weak solution of \eqref{2}, taking in \eqref{19} test functions $\psi\equiv \phi\in H^1$ for an arbitrary $T>0$ we get \eqref{18}, concluding the proof.
\end{proof}
Equation \eqref{18}  has a variational nature, as the next result highlights.
\begin{lem}\label{Lemma 2} Il \eqref{3} holds the functional $I$ defined in \eqref{6} belongs to $C^1(H^1)$ and its critical points coincide with the weak solutions of \eqref{1}.
\end{lem}
\begin{proof}By classical arguments, see \cite[Chapter 1, Theorem~2.9, p. 22]{ambrosettiprodi}, the potential operator $F:H^1\to\R$, defined by $F(u)=\frac 1p \|u\|_{p,\Gamma_1}^p$, is Fr\'{e}chet differentiable and, for all $u,\phi\in H^1$ one has
\begin{equation}\label{20}
  \langle F'(u),\phi\rangle_{H^1}=\int_{\Gamma_1}|u|^{p-2}u\phi.
\end{equation}
Consequently, since
\begin{equation}\label{20bis}
 I(u)=\tfrac 12\|u\|_{H^1}^2-F(u)\qquad\text{for all $u\in H^1$,}
\end{equation}
trivially $I\in C^1(H^1)$ and
\begin{equation}\label{21}
\langle I'(u),\phi\rangle_{H^1}=(u,\phi)_{H^1}-\langle F'(u),\phi\rangle_{H^1}\quad\text{for all $u,\phi\in H^1$.}
\end{equation}
By \eqref{10}, \eqref{20} and \eqref{21} one immediately  gets that \eqref{18} is rewritten as $I'(u)=0$.
\end{proof}
We now establish some geometrical properties of the functional $I$.
\begin{lem}\label{Lemma 3}Il \eqref{3} holds the functional $I$ satisfies the assumptions i)--iii) of Theorem~\ref{Theorem 3} and
iv) of Theorem~\ref{Theorem 4}.
\end{lem}
\begin{proof}
By \eqref{6} trivially $I(0)=0$, proving i). To prove ii) and iv) we remark that, by Sobolev Embedding Theorem there is $c_4=c_4(p,\Omega)>0$ such that
$$\|u\|_{p,\Gamma}\le c_4 \|u\|_{H^1(\Gamma)}\qquad\text{for all $u\in H^1(\Gamma)$.}$$
Consequently, by Lemma~\ref{lemma 0}, there is $c_5=c_5(p,\Omega,\Gamma_0)>0$ such that
\begin{equation}\label{21bis}
  \|u\|_{p,\Gamma_1}\le c_5 \|u\|_{H^1}\qquad\text{for all $u\in H^1$.}
\end{equation}
Hence, by \eqref{6}, for all $u\in H^1$ we have
$$I(u)=\frac 12 \|u\|_{H^1}^2-\frac 1p \|u\|_{p,\Gamma_1}^p\ge
\left(\frac 12-\frac{c_5^p}p \|u\|_{H^1}^{p-2}\right)\|u\|_{H^1}^2,$$
from which ii) and iv) follow, by taking for example
$$\rho=\left(\frac p{4c_5^p}\right)^{\frac 1{p-2}},\quad\text{and}\quad \alpha=\frac 14 \left(\frac p{4c_5^p}\right)^{\frac 2{p-2}}.$$
To prove iii) we remark that, for any $u\in H^1$ such that $u_{|\Gamma_1}\not=0$ and $s>0$ we have
$$I(su)=\tfrac 12 \|u\|_{H^1}^2s^2-\tfrac 1p \|u\|_{p,\Gamma_1}^p s^p\to-\infty\quad\text{as $s\to\infty$.}$$
\end{proof}
The last relevant property of $I$ is given by the following result.
\begin{lem}\label{Lemma 4} Il \eqref{3} holds the functional $I$ satisfies the (PS) condition.
\end{lem}
\begin{proof}Let $(u_n)_n$ be a (PS) sequence in $H^1$. Then there are $c_6,c_7\ge 0$, depending on  $(u_n)_n$, such that
\begin{equation}\label{22}
 I(u_n)\le c_6,\qquad\text{and}\quad |\langle I'(u_n),u_n\rangle_{H^1}|\le c_7 \|u_n\|_{H^1}\quad\text{for all $n\in\N$.}
\end{equation}
Since, by \eqref{6} and \eqref{21},
$$pI(u_n)-\langle I'(u_n),u_n\rangle_{H^1}=\left(\tfrac p2-1\right)\|u_n\|_{H^1}^2,$$
by \eqref{22} we get
$$\left(\tfrac p2-1\right)\|u_n\|_{H^1}^2\le p c_6+c_7\|u_n\|_{H^1},$$
from which one immediately yields that $(u_n)_n$ is bounded in $H^1$.
Consequently, up to a subsequence, $u_n\rightharpoonup u$ in $H^1$. To prove that the convergence is actually strong we remark that, by \eqref{21}, for all $\phi\in H^1$ we have
$$(u_n-u,\phi)_{H^1}=\langle I'(u_n)-I'(u),\phi\rangle_{H^1}+\int_{\Gamma_1}(|u_n|^{p-2}u_n-|u|^{p-2}u)\phi
$$
so taking $\phi=u_n-u$ we get
\begin{multline}\label{23}
\|u_n-u\|_{H^1}^2= \langle I'(u_n),u_n-u\rangle_{H^1}-\langle I'(u),u_n-u\rangle_{H^1}\\
 +\int_{\Gamma_1}(|u_n|^{p-2}u_n-|u|^{p-2}u)(u_n-u).
\end{multline}
The first two terms in the right hand side  of \eqref{23} converges to $0$  since $I'(u_n)\to 0$ in $(H^1)'$ and  $u_n-u\rightharpoonup 0$, hence being norm bounded. As to the third  term in it, since the embedding $H^1(\Gamma)\hookrightarrow L^p(\Gamma)$ is compact and the operator $u\mapsto u_{|\Gamma}$ from $H^1$ to $H^1(\Gamma)$ is bounded, up to subsequence we have ${u_n}_{|\Gamma}\to u_{|\Gamma}$  strongly in $L^p(\Gamma)$. By standard properties of the Nemitskii operators, see \cite[Chapter 1, Theorem~2.2, p. 16]{ambrosettiprodi}, we also have $|{u_n}_{|\Gamma}|^{p-2}{u_n}_{|\Gamma}\to |u_{|\Gamma}|^{p-2}u_{|\Gamma}$  strongly in $L^{p'}(\Gamma)$. Consequently the third term in the right hand side  of \eqref{23} converges to $0$ goes to $0$ by H\"{o}lder inequality.
\end{proof}
We can then give the proof of our first main result.
\begin{proof}[\bf Proof of Theorem~\ref{Theorem 1}] By simply combining Lemmas~\ref{Lemma 1}--\ref{Lemma 4},
Theorem~\ref{Theorem 4}, with $I=\mathcal{I}$, and the fact that $I$ is even we get the statement, but for the critical level, since by Theorem~\ref{Theorem 4} we have that $I(u)=I(-u)=c$, where $c$ is given by \eqref{14}. To complete the proof we then have to recognize that $c=d_1=d_2$, where
$$d_1=\inf_{u\in H^1, u_{|\Gamma}\not\equiv 0}\sup_{\lambda>0} I(\lambda u),\qquad\text{and}\qquad d_2\inf_{u\in H^1\setminus\{0\}}\sup_{\lambda>0} I(\lambda u).$$
Since, for any $\lambda>0$,
\begin{equation}\label{24}
I(\lambda u)=\tfrac 12 \|u\|_{H^1}^2\lambda^2-\tfrac 1p \|u\|_{p,\Gamma_1}^p \lambda^p,
\end{equation}
trivially $\sup\limits_{\lambda>0} I(\lambda u)<\infty$ when $u_{|\Gamma}\not\equiv 0$, while $\sup\limits_{\lambda>0} I(\lambda u)=\infty$ when $u\in H^1\setminus\{0\}$ and $u_{|\Gamma}\equiv 0$. Hence $d_1=d_2=d$, and consequently it remains to prove that $d=c$.

This fact is well--known for similar problems, see for example \cite[Chapter 8, p. 117]{ambrosettimalchiodi} and  \cite{bellazzinivisciglia,marcosmedeiros,jeanjeantanaka2, jeanjeantanaka1,rendiconti}. Here we shall essentially conveniently adapt the argument in \cite{rendiconti}.
By \eqref{24}, for any $u\in H^1$ with $u_{|\Gamma}\not\equiv 0$, the function $\lambda\mapsto I(\lambda u)$ has a unique critical point $\lambda_u>0$, $\max\limits_{\lambda>0} I(\lambda u)=I(\lambda_u u)$ and $I(\lambda u)\to-\infty$ as $\lambda\to\infty$, for any such $u$, defining $\sigma_u\in C([0,1];H^1)$ by $\sigma_u(s)=Rsu$ for $s\in [0,1]$, with $R>0$ so large that $I(Ru)<0$, we have $\sigma_u\in \Sigma$ and consequently
$I(\lambda_u u)=\max\limits_{t\in [0,1]}I(\sigma_u(t))\ge c$.  Hence $d\ge c$. On the other hand, if $u$ is a critical point of $I$ with $I(u)=c$,  already found above, by \eqref{21} we have $\|u\|_{H^1}^2=\|u\|_{p,\Gamma_1}^p$. Then, since $u\not=0$, we also get that $u_{|\Gamma}\not\equiv 0$. Consequently, since
$\frac d{d\lambda}I(\lambda u)=\langle I'(u),u\rangle_{H^1}$, we have $\lambda_u=1$ and consequently
$c=I(u)=\max\limits_{\lambda>0} I(\lambda u)\ge d$, completing the proof.
\end{proof}
We now turn to prove Theorem~\ref{Theorem 1.1}. We remark at first that the number $B$ defined in \eqref{7bis} is finite because of the estimate \eqref{21bis}. We now introduce the auxiliary functional $K\in C^1(H^1)$ given by
\begin{equation}\label{42}
 K(u)=\langle I'(u),u\rangle_{H^1}=\|u\|_{H^1}^2-\|u\|_{p,\Gamma_1}^p.
\end{equation}
The key point in the proof of Theorem~\ref{Theorem 1.1} is the following result, of possible independent interest.
\begin{lem}\label{lemma 9} Let \eqref{3} hold and $\lambda_1$, $\lambda_2$ be given by \eqref{7ter}. Then \eqref{43} holds.
Moreover, for any $u\in H^1$ such that $u_{|\Gamma}\not\equiv 0$ and $I(u)\le d$, the following implications hold:
\begin{equation}\label{44}
\begin{alignedat}3
&K(u)\ge 0\quad & \Longleftrightarrow & \quad \|u\|_{H^1}\le\lambda_1\quad & \Longleftrightarrow& \quad \|u\|_{p,\Gamma_1}\le \lambda_2, \\
&K(u)\le 0\quad & \Longleftrightarrow & \quad \|u\|_{H^1}\ge\lambda_1\quad & \Longleftrightarrow& \quad \|u\|_{p,\Gamma_1}\ge \lambda_2.
\end{alignedat}
\end{equation}
\end{lem}
\begin{proof} To prove \eqref{43} we remark that, for any $u\in H^1$ with $u_{|\Gamma}\not\equiv 0$ an easy calculations shows that
$$\max_{\lambda>0}I(\lambda u)=I(\lambda_u u)\quad\text{where}\quad \lambda_u=\|u\|_{H^1}^{2/(p-2)}\|u\|_{p,\Gamma_1}^{-p/(p-2)}$$
so that
$$\max_{\lambda>0}I(\lambda u)=\left(\frac 12 -\frac 1p\right)\left(\frac {\|u\|_{H^1}}{\|u\|_{p,\Gamma_1}}\right)^{\frac {2p}{p-2}}$$
and consequently
$$d=\left(\frac 12 -\frac 1p\right)\left(\sup_{u\in H^1, u_{|\Gamma\not\equiv 0}}\frac {\|u\|_{H^1}}{\|u\|_{p,\Gamma_1}}\right)^{\frac {-2p}{p-2}}.$$
By using \eqref{7bis} and \eqref{7ter} in the last formula we get \eqref{43}.

Now let $u\in H^1$ such that $u_{|\Gamma}\not\equiv 0$ and $I(u)\le d$. To prove \eqref{44} we shall first prove the implications
\begin{equation}\label{45}
K(u)\ge 0\quad  \Longrightarrow  \quad \|u\|_{H^1}\le\lambda_1\quad  \Longrightarrow \quad \|u\|_{p,\Gamma_1}\le \lambda_2 \quad  \Longrightarrow \quad K(u)\ge 0.
\end{equation}
If $K(u)\ge 0$, supposing by contradiction that $\|u\|_{H^1}>\lambda_1$, by \eqref{20bis} we get
$$I(u)\ge \left(\tfrac 12 -\tfrac 1p\right)\|u\|_{H^1}^2>\left(\tfrac 12 -\tfrac 1p\right)\lambda_2^2=d,$$
a contradiction.
If $\|u\|_{H^1}\le \lambda_1$, since by \eqref{7bis} we have
\begin{equation}\label{46}
  \|u\|_{p,\Gamma_1}\le B \|u\|_{H^1}\qquad\text{for all $u\in H^1$,}
\end{equation}
and $\lambda_2=B\lambda_1$, we immediately get $\|u\|_{p,\Gamma_1}\le \lambda_2$.
If $\|u\|_{p,\Gamma_1}\le \lambda_2$, by \eqref{7ter} and \eqref{46} we have
$$\|u\|_{p,\Gamma_1}^p\le \lambda_2^{p-2}\|u\|_{p,\Gamma_1}^2=B^{-2}\|u\|_{p,\Gamma_1}^2\le \|u\|_{H^1}^2,$$
so $K(u)\ge 0$, concluding the proof of \eqref{45}.

To complete the  proof of \eqref{44} we are then going to prove the further implications
\begin{equation}\label{45}
K(u)\le 0\quad  \Longrightarrow  \quad \|u\|_{p,\Gamma_1}\ge \lambda_2\quad  \Longrightarrow \quad \|u\|_{H^1}\ge\lambda_1 \quad  \Longrightarrow \quad K(u)\le 0.
\end{equation}
If $K(u)\le 0$, by \eqref{46}, we have
$$B^{-2}\|u\|_{p,\Gamma_1}^2\le \|u\|_{H^1}^2\le \|u\|_{p,\Gamma_1}^p$$
which, since $u_{|\Gamma}\not\equiv 0$, gives $\|u\|_{p,\Gamma_1}\ge B^{-2/(p-2)}=\lambda_2$.
If $\|u\|_{p,\Gamma_1}\ge \lambda_2$  by \eqref{46} we have $B\|u\|_{H^1}\ge \lambda_1$, i.e. $\|u\|_{H^1}\ge B^{-1}\lambda_2=\lambda_1$.
If $\|u\|_{H^1}\ge\lambda_1$ then, assuming by contradiction that $K(u)>0$, i.e. $\|u\|_{p,\Gamma_1}^p<\|u\|_{H^1}^2$, by \eqref{20bis} we get
$$I(u)=\tfrac 12  \|u\|_{H^1}^2-\tfrac 1p \|u\|_{p,\Gamma_1}^p>\left(\tfrac 12 -\tfrac 1p\right)\|u\|_{H^1}^2
\ge \left(\tfrac 12 -\tfrac 1p\right)\lambda_1^2=d,$$
a contradiction, so $K(u)\le 0$ and the proof is complete.
\end{proof}
We can now prove our second main result.

\begin{proof}[\bf Proof of Theorem~\ref{Theorem 1.1}] Formula \eqref{43} has been already proved in Lemma~\ref{lemma 9}. Moreover, if $u$ is a weak solution of \eqref{1} with $I(u)=d$, by Lemma~\ref{Lemma 2} and \eqref{42} we have $K(u)=0$, and since $d>0$ necessarily
$u\not\equiv 0$. Hence $\|u\|_{p,\Gamma_1}^p=\|u\|_{H^1}^2>0$, so by Lemma~\ref{lemma 9} we immediately get \eqref{7quater}.

Now let $u$ a non--trivial weak solution of \eqref{1}. By Lemma~\ref{Lemma 2} and \eqref{42} we have $K(u)=0$. Then, by \eqref{20bis} and \eqref{46} we get
\begin{equation}\label{48}
 I(u)\ge \frac 12 \|u\|_{H^1}^2-\frac {B^p}p\|u\|_{H^1}^p
\end{equation}
and, since $K(u)=0$, also
\begin{equation}\label{49}
 I(u)=\tfrac 12 \|u\|_{H^1}^2-\tfrac 1p\|u\|_{H^1}^2=\left(\tfrac 12 -\tfrac 1p\right)\|u\|_{H^1}^2.
\end{equation}
By combining \eqref{48} and \eqref{49} we get $B^p\|u\|_{H^1}^p\ge \|u\|_{H^1}^2$, so since $u\not\equiv 0$ we have $\|u\|_{H^1}\ge B^{-p/(p-2)}=\lambda_1$. Then, using \eqref{49} again we obtain $I(u)\ge \left(\tfrac 12 -\tfrac 1p\right)\lambda_1^2=d$, concluding the proof.
\end{proof}
\section{Multiplicity.} \label{section 4}
The aim of this section is to prove Theorem~\ref{Theorem 2} and the strategy of the proof consists in applying Theorem~\ref{Theorem 5}. Unfortunately the simple variational setting used in the previous section is not adequate for this purpose. Indeed the functional $I$ defined in \eqref{6} does not satisfy assumption (v) in Theorem~\ref{Theorem 5}., as it is evident by considering its restriction to the space $H^1_0(\Omega)$.

To introduce a convenient setting we first recall that. by standard elliptic theory, see for example \cite[Chapter I]{GiraultRaviart}, for any $v\in H^{1/2}(\Gamma)$ the nonhomogeneous Dirichlet problem
\begin{equation}\label{26}
\begin{cases}
\Delta u=0&\quad\text{in $\Omega$,}\\
u=v &\quad\text{on $\Gamma$,}
\end{cases}
\end{equation}
has a unique solution $u\in H^1(\Omega)$, continuously depending on $v$ in the topologies of the respective spaces, where the Laplace equation $\Delta u=0$ is taken in the sense of distributions, or equivalently in the space $H^{-1}(\Omega)$, that is
\begin{equation}\label{27}
  \int_\Omega \nabla u\nabla \phi=0\qquad\text{for all $\phi\in H^1_0(\Omega)$,}
\end{equation}
while the boundary condition is taken in the trace sense. Trivially $u\in H^1$ when $v\in H^1_{\Gamma_0}(\Gamma)$, see \eqref{++},
so obtaining the bounded linear Dirichlet operator $v\mapsto u$ from $H^1_{\Gamma_0}(\Gamma)$ into $H^1$.
Trivially its range is the closed subspace $\mathbb{A}^1$ of $H^1$ defined by
\begin{equation}\label{28}
\mathbb{A}^1=\{u\in H^1: \quad\text{\eqref{27} holds}\}.
\end{equation}
Hence, denoting $u=\mathbb{D}v$, we get the bijective isomorphism
\begin{equation}\label{29}
\mathbb{D}\in\mathcal{L}(H^1_{\Gamma_0}(\Gamma);\mathbb{A}^1),\quad\text{with $\mathbb{D}^{-1}=\Tr_{|\mathbb{A}^1}$.}
\end{equation}
The starting point of the analysis is that the space $\mathbb{A}^1$ is a natural constraint for problem \eqref{1}, since equations \eqref{1}$_1$ and \eqref{1}$_2$ hold in it. To write equation \eqref{1}$_3$ in a weak form we shall use the realization $-\Delta_{\Gamma_{1(D)}}$ of the Laplace--Beltrami operator introduced in \eqref{++} as well as the Dirichlet--to--Neumann operator $\mathcal{A}\in\mathcal{L}(H^1_{\Gamma_0}(\Gamma);[H^1_{\Gamma_0}(\Gamma)]' )$ given by
\begin{equation}\label{31}
\langle \mathcal{A}v,w\rangle_{H^1_{\Gamma_0}(\Gamma)}=\int_\Omega \nabla(\mathbb{D}v)\nabla(\mathbb{D}w)
\quad\text{for all $v,w\in H^1_{\Gamma_0}(\Gamma)$.}
\end{equation}
By \eqref{28} and \eqref{29}, since problems \eqref{26} has a unique solution, we also get that
\begin{equation}\label{32}
\langle \mathcal{A}v,w\rangle_{H^1_{\Gamma_0}(\Gamma)}=\int_\Omega \nabla(\mathbb{D}v)\nabla\phi
\end{equation}
for all $v,w\in H^1_{\Gamma_0}(\Gamma)$ and $\phi\in H^1(\Omega)$ such that $\phi_{|\Gamma}=w$.
One can then abstractly write equation \eqref{1}$_3$ in the space $[H^1_{\Gamma_0}(\Gamma)]' $ as follows.
\begin{definition}\label{Definition 4}Let \eqref{3} hold. We say that $v\in H^1_{\Gamma_0}(\Gamma)$ is a solution of
\begin{equation}\label{33}
 -\Delta_{\Gamma_{1(D)}}v+\mathcal{A}v=|v|^{p-2}v\qquad\text{in $[H^1_{\Gamma_0}(\Gamma)]'$}
\end{equation}
if
\begin{equation}\label{34}
\int_\Omega \nabla(\mathbb{D}v)\nabla(\mathbb{D}\psi)+\int_{\Gamma_1}(\nabla_\Gamma v,\nabla_\Gamma \psi)_\Gamma-\int_{\Gamma_1}|v|^{p-2}v\psi=0
\end{equation}
for all $\psi\in H^1_{\Gamma_0}(\Gamma)$.
\end{definition}
We immediately get the following result.
\begin{lem}\label{Lemma 5}
Let $u\in H^1$. Then $u$ is a weak solution of \eqref{1} if and only if $u\in\mathbb{A}^1$ and $v=u_{|\Gamma}$ is a solution of \eqref{33}. Moreover in this case $u=\mathbb{D}v$.
\end{lem}
\begin{proof}
If $u\in H^1$ is a weak solution of \eqref{1}, by \eqref{18} one immediately gets that $u\in\mathbb{A}^1$, so $v=u_{|\Gamma}$ satisfies \eqref{34} and $u=\mathbb{D}v$. Conversely, if  $v\in H^1_{\Gamma_0}(\Gamma)$ satisfies \eqref{34} then $u:=\mathbb{D}v\in \mathbb{A}^1$ and, by combining \eqref{31}, \eqref{32} and \eqref{34} we get \eqref{18}.
\end{proof}
Trivially equation \eqref{34} has a variational structure. To write it down  we introduce the functional $J:H^1_{\Gamma_0}(\Gamma)\to\R$ given by
\begin{equation}\label{35}
 J(v)=\frac 12 \int_\Omega |\nabla(\mathbb{D}v)|^2+\frac 12 \int_{\Gamma_1} |\nabla_\Gamma v|_\Gamma^2-\frac 1p \int_{\Gamma_1}|v|^p.
\end{equation}
The following result points out some trivial properties of $J$.
\begin{lem}\label{Lemma 6} Il \eqref{3} holds we have $J\in C^1(H^1_{\Gamma_0}(\Gamma))$, its Fr\'{e}chet derivative being given
for all $v,\psi\in H^1_{\Gamma_0}(\Gamma)$ by
\begin{equation}\label{36}
 \langle J'(v),\psi\rangle_{H^1_{\Gamma_0}(\Gamma)}=\int_\Omega \nabla(\mathbb{D}v)\nabla(\mathbb{D}\psi)+\int_{\Gamma_1}(\nabla_\Gamma v,\nabla_\Gamma \psi)_\Gamma-\int_{\Gamma_1}|v|^{p-2}v\psi.
\end{equation}
Consequently solutions of \eqref{33} are exactly critical points of $J$. Moreover $J=I\cdot \mathbb{D}$, $J$ is even and $J(0)=0$.
\end{lem}
\begin{proof}
Trivially $J$ is even, $J(0)=0$ and $J=I\cdot \mathbb{D}$. Hence, since $I\in C^1(H^1)$, we also have  $J\in C^1(H^1_{\Gamma_0}(\Gamma))$ and \eqref{36} holds true. The correspondence between solutions of \eqref{33} and critical points of $J$
immediately follows by \eqref{34} and \eqref{36}.
\end{proof}
The next result shows that $J$ satisfies the  geometrical assumptions of Theorem~\ref{Theorem 5} which are still to be checked.
\begin{lem}\label{Lemma 7} Il \eqref{3} holds the functional $J$ satisfies assumption ii) in Theorem~\ref{Theorem 3} and assumption v) in Theorem~\ref{Theorem 5}.
\end{lem}
\begin{proof} Since $J=I\cdot \mathbb{D}$, by \eqref{21bis} there is $c_8=c_8(p,\Omega,\Gamma_0)>0$ such that
$$\|v\|_{p,\Gamma_1}\le c_8 \|\mathbb{D}v\|_{H^1}\qquad\text{for all $v\in H^1_{\Gamma_0}(\Gamma)$.}$$
Consequently, by \eqref{11} and \eqref{35}, we get the estimate
$$J(v)\ge
\left(\frac 12-\frac{c_8^p}p \|\mathbb{D}v\|_{H^1}^{p-2}\right)\|\mathbb{D}v\|_{H^1}^2\qquad\text{for all $v\in H^1_{\Gamma_0}(\Gamma)$,}$$
and then
\begin{equation}\label{38}
  J(v)\ge \frac 14 \|\mathbb{D}v\|_{H^1}^2\quad\text{for all $v\in H^1_{\Gamma_0}(\Gamma)$ such that}
  \quad \|\mathbb{D}v\|_{H^1}\le \left(\frac p{4c_8^p}\right)^{\frac 1{p-2}}.
\end{equation}
Since, by \eqref{29}, $\|\mathbb{D}(\cdot)\|_{H^1}$ is on $H^1_{\Gamma_0}(\Gamma)$ a norm equivalent to the one inherited from $H^1(\Gamma)$, assumption ii) in Theorem~\ref{Theorem 3}, for suitable $\rho,\alpha>0$, trivially follows from \eqref{38}.

To prove that also assumption v) in Theorem~\ref{Theorem 5} holds true let $Y$ be any finite dimensional subspace of $H^1_{\Gamma_0}(\Gamma)$. Since on $Y$ all norms are equivalent there are $c_9=c_9(p,Y)>0$ and  $c_{10}=c_{10}(p,Y)>0$ such that
\begin{equation}\label{38bis}
c_9\|v\|_{p,\Gamma_1}\le \|v\|_{H^1(\Gamma)}\le c_{10} \|v\|_{p,\Gamma_1}\quad\text{for all $v\in Y$. }
\end{equation}
Consequently, using again \eqref{29}, there is $c_{11}=c_{11}(p,Y,\Omega,\Gamma_0)>0$  such that
$$\|\mathbb{D} v\|_{H^1}\le c_{11}\|v\|_{p,\Gamma_1}\qquad\text{for all $v\in Y$.}$$
Using it in \eqref{35} we then get
$$J(v)=\frac 12 \|\mathbb{D}v\|_{H^1}^2-\frac 1p \|v\|_{p,\Gamma_1}^p\le
\frac 1{2c_{11}^2} \|v\|_{p,\Gamma_1}^2-\frac 1p \|v\|_{p,\Gamma_1}^p,
$$
from which, since $p>2$, one gets that  $J(v)\le 0$ for $\|v\|_{p,\Gamma_1}$ sufficiently large.
By \eqref{38bis} then assumption v) in Theorem~\ref{Theorem 5} holds true.
\end{proof}
To complete checking  of the assumptions of Theorem~\ref{Theorem 5} we give the following result.
\begin{lem}\label{Lemma 8}  Il \eqref{3} holds the functional $J$ satisfies the (PS) condition.
\end{lem}
\begin{proof}We shall prove the statement by using  Lemma~\ref{Lemma 4}. With this aim we make some preliminary remarks concerning the space $H^1$ and the functional $I$.
At first $H^1$ admits the orthogonal (with respect to $(\cdot,\cdot)_{H^1}$ given in \eqref{10}) splitting
\begin{equation}\label{40}
  H^1=\mathbb{A}^1\oplus H^1_0(\Omega),
\end{equation}
the respective orthogonal projectors $\Pi_{\mathbb{A}^1}\in\mathcal{L}(H^1,\mathbb{A}^1)$ and $\Pi_{H^1_0(\Omega)}\in\mathcal{L}(H^1,H^1_0(\Omega))$ being given by
$$\Pi_{\mathbb{A}^1}u=\mathbb{D}u_{|\Gamma}\quad\text{and}\quad \Pi_{H^1_0(\Omega)}u=u-\mathbb{D}u_{|\Gamma}\quad\text{for all $u\in H^1$.} $$
Using this splitting we can rewrite \eqref{21}, for any $u\in \mathbb{A}^1$,  as
\begin{equation}\label{41}
\begin{aligned}
 \langle I'(u),\phi+\psi\rangle_{H^1}=&(u,\phi+\psi)_{H^1}+\int_{\Gamma_1}|u|^{p-2}u(\phi+\psi)\\
 =&(u,\phi)_{H^1}+\int_{\Gamma_1}|u|^{p-2}u\phi\\
 =&\langle I'_{\mathbb{A}^1}(u),\phi\rangle_{\mathbb{A}^1}\quad\text{for all $\phi\in\mathbb{A}^1$ and $\psi\in H^1_0(\Omega)$,}
\end{aligned}
\end{equation}
where $I_{\mathbb{A}^1}=I_{|\mathbb{A}^1}$, that is $I_{\mathbb{A}^1}$ is the restriction of $I$ to $\mathbb{A}^1$.
\begin{footnote}{once one recognizes that critical points of $I$ belong to $\mathbb{A}^1$,  formula \eqref{41} gives a different  proof of Lemma~\ref{Lemma 5} }\end{footnote}

To prove the statement let now $(v_n)_n$  be a (PS) sequence for $J$.
Since $J=I\cdot\mathbb{D} =I_{\mathbb{A}^1}\cdot \mathbb{D}$, by \eqref{29} we get that $(\mathbb{D}v_n)_n$ is a (PS) sequence for
the functional $I_{\mathbb{A}^1}$ in $\mathbb{A}^1$. Hence $(I(\mathbb{D}v_n))_n$ is bounded in $H^1$ and, using \eqref{41},
$I'(\mathbb{D}v_n)\to 0$ in $[H^1]'$, so $(\mathbb{D}v_n)_n$ is a (PS) sequence for $I$.

By  Lemma~\ref{Lemma 4} then, up to a subsequence, $(\mathbb{D}(v_n))_n$ converges in $H^1$ which, using  \eqref{29} again, concludes the proof.
\end{proof}
We can finally give the proof of our last main result.
\begin{proof}[\bf Proof of Theorem~\ref{Theorem 2}]
By Lemmas~\ref{Lemma 6}--\ref{Lemma 8} the functional $J$ satisfies the assumptions of Theorem~\ref{Theorem 5}, so it possesses
a sequence $(v_n)_n$ of critical points such that $J(v_n)\to\infty$. By Lemma~\ref{Lemma 5} the sequence $(u_n)_n$ given by $u_n=\mathbb{D}v_n$ for all $n\in\N$ is thus a sequence of critical points of $I$ with $I(u_n)=J(v_n)\to\infty$. By Lemma~\ref{Lemma 3} the proof is thus concluded.
\end{proof}

\def\cprime{$'$}
\providecommand{\href}[2]{#2}
\providecommand{\arxiv}[1]{\href{http://arxiv.org/abs/#1}{arXiv:#1}}
\providecommand{\url}[1]{\texttt{#1}}
\providecommand{\urlprefix}{}

\end{document}